\documentclass[3p,times]{elsarticle}
\usepackage{amsmath}
\usepackage{mathrsfs}
\usepackage{enumerate}
\usepackage[colorlinks=true]{hyperref}
\newtheorem{thm}{Theorem}
\newtheorem{lem}{Lemma}
\newtheorem{coy}{Corollary}
\newtheorem{prn}{Proposition}

\newdefinition{rmk}{Remark}
\newdefinition{exe}{Example}
\newdefinition{qun}{Question}

\newcommand{\pf}{\noindent \textbf{Proof.}\quad}
\newcommand{\epf}{\hspace{\stretch{1}}$\blacksquare$}

\begin{document}

\begin{frontmatter}
\title{Topological dynamics of Zadeh's extension on the space of upper semi-continuous
fuzzy sets\tnoteref{label1}}
\tnotetext[label1]{Project supported by
the scientific research starting project
of Southwest Petroleum University (No. 2015QHZ029), the Scientific Research
Fund of the Sichuan Provincial Education Department (No. 14ZB0007) and NSFC (No. 11401495).}
\author[rvt]{Xinxing Wu}
\ead{wuxinxing5201314@163.com}

\address[rvt]{School of Sciences and Research Center for
Mathematics and Mechanics,
Southwest Petroleum University, \\
Chengdu, Sichuan, 610500, P.R. China
}

\author[rvt2]{Xiong Wang\corref{cor1}}
\ead{wangxiong8686@szu.edu.cn}
\cortext[cor1]{Corresponding author}

\address[rvt2]{Institute for Advanced Study, Shenzhen University,\\
 Nanshan District Shenzhen, Guangdong, China 
}


\begin{abstract}
In this paper, some characterizations about transitivity, mildly mixing property, $\mathbf{a}$-transitivity,
equicontinuity, uniform rigidity and proximality of Zadeh's extensions restricted on some invariant closed subsets of the space of all upper semi-continuous fuzzy sets with the level-wise metric are obtained. In particular, it is proved that a dynamical system is weakly mixing (resp., mildly mixing, weakly mixing and $\mathbf{a}$-transitive, equicontinuous, uniformly rigid) if and only if
the Zadeh's extension is transitive (resp., mildly mixing, $\mathbf{a}$-transitive, equicontinuous, uniformly rigid).
\end{abstract}
\begin{keyword}
Zadeh's extension, transitivity, mildly mixing, equicontinuous, rigidity.
\MSC[2010] 03E72, 37A25, 54A40, 54C60, 54H20.

\end{keyword}

\end{frontmatter}


\section{Introduction}

A {\it dynamical system} is a pair $(X, T)$,
where $X$ is a nontrivial compact metric space with a metric $d$ and
$T: X\longrightarrow X$ is a continuous surjection. A nonempty invariant closed subset $Y
\subset X$ (i.e., $T(Y) \subset Y$) defines naturally a {\it subsystem} $(Y,
T|_{Y})$ of $(X, T)$. Throughout this paper, let $\mathbb{N}=\{1, 2, 3, \ldots\}$,
$\mathbb{Z}^{+}=\{0, 1, 2, \ldots\}$ and $I=[0, 1]$. For any $n\in \mathbb{N}$,
write $(X^{n}, T^{(n)})$ as the {\it $n$-fold product system}
$(\underbrace{X\times\dots\times X}\limits_{n}, \underbrace{T\times\dots\times T}\limits_{n})$.

\medskip

Sharkovsky's amazing discovery \cite{Sharkovsky1964}, as well
as Li and Yorke's famous work which introduced the concept of
`chaos' known as Li-Yorke chaos today
\cite{LY75}, have provoked the recent rapid advancement of research on
discrete chaos theory. The essence of Li-Yorke chaos is the
existence of uncountable scrambled sets. Another well-known definition
of chaos was given by Devaney \cite{D89}, according to which
a continuous map $T$ is said to be {\it chaotic in the sense of
Devaney} if it satisfies the following three properties:
\begin{enumerate}[(1)]
\item $T$ is {\it (topologically) transitive}, i.e., for every pair of
nonempty open sets $U, V$ of $X$, there exists
$n\in \mathbb{Z}^{+}$ such that
$T^{n}(U)\cap V\neq \emptyset$;

\item The set of periodic points of $T$ is dense in $X$;

\item\label{D-3} $T$ has {\it sensitive dependence on initial conditions}
(briefly, {\it is sensitive}), i.e., there exists $\varepsilon >0$
such that for any $x\in X$ and any neighborhood $U$ of $x$,
there exist $y\in U$ and $n\in \mathbb{Z}^{+}$ satisfying
$d(T^{n}(x), T^{n}(y))>\varepsilon$.
\end{enumerate}

Banks et al. \cite{Banks1992} proved that every transitive map whose periodic points are dense
in $X$ is sensitive, which implies that the above condition (\ref{D-3}) is redundant, while Huang and Ye
\cite{Huang-Ye2002} showed that every transitive map containing a periodic point is chaotic
in the sense of Li-Yorke.

Given a dynamical system $(X, T)$, one can naturally obtain three associated systems induced by $(X, T)$.
The first one is $({K}(X), T_{K})$ on the hyperspace ${K}(X)$ consisting
of all nonempty closed subsets of $X$ with the Hausdorff metric. The second one
is $(M(X), T_M)$ on the space $M(X)$ consisting of all Borel probability measures with
the Prohorov metric. And the last one is its Zadeh's extension $(\mathbb{F}_{0}(X), T_{F})$
(more generally $g$-fuzzification $(\mathbb{F}_{0}(X), T_{F}^{g})$) on the space $\mathbb{F}_{0}(X)$
consisting of all nonempty upper semi-continuous fuzzy sets with the level-wise metric
induced by the extended Hausdorff metric.
A systematic study on the connections between dynamical properties of
$(X, T)$ and its two induced systems $(K(X), T_{K})$ and
$(M(X), T_{M})$ was initiated by Bauer and Sigmund in \cite{BS75}, and later has been widely
developed by several authors. For more results on this topic,
one is referred to \cite{Banks2005,Illanes1999,Gu2007,Li2014,LOYZ,Liao2006,Liu,SN10,WOC,WWC09,WWC,WC}
and references therein.

\medskip

Rom\'{a}n-Flores and  Chalco-Cano \cite{Roman2008} studied
some chaotic properties (for example,
transitivity, sensitive dependence, periodic density) for the Zadeh's extension
of a dynamical system. Then, Kupka \cite{Kupka2011} investigated the relations
between Devaney chaos in the original system and in the Zadeh's extension
and proved that Zadeh's extension is periodically dense in $\mathbb{F}(X)$
(resp. $\mathbb{F}^{\geq\lambda}(X)$ for any $\lambda\in (0, 1]$) if and only if so is
$(K(X), T_{K})$ (see Lemma \ref{piecewise}). Recently, Kupka \cite{Kupka2011-1} introduced
the notion of $g$-fuzzification to generalize Zadeh's extension and obtained some
basic properties of $g$-fuzzification. In \cite{Kupka2014}, he continued in studying chaotic properties
(for example, Li-Yorke chaos, distributional chaos, $\omega$-chaos, transitivity,
total transitivity, exactness, sensitive dependence, weakly mixing, mildly mixing, topologically mixing)
of $g$-fuzzification and showed that if the $g$-fuzzification
$(\mathbb{F}^{=1}(X), T^{g}_{F})$ has the property $P$, then $(X, T)$
also has the property $P$, where $P$ denotes the following properties: exactness, sensitive dependence,
weakly mixing, mildly mixing, or topologically mixing. Meanwhile, he posed the following question:

\begin{qun}\cite{Kupka2014}\label{Question 1}
Does the $P$-property of $(X, T)$ imply the
$P$-property of $(\mathbb{F}^{=1}(X), T^{g}_{F})$?
\end{qun}

We \cite{Wu} obtained a sufficient condition on $g\in D_{m}(I)$\footnote{$D_{m}(I)$ is the set of all nondecreasing right-continuous functions
$g: I\longrightarrow I$ with $g(0)=0$ and $g(1)=1$.} to ensure that for every dynamical system
$(X, T)$, its $g$-fuzzification $(\mathbb{F}^{=1}(X), T_{F}^{g})$ is not transitive
(thus, not weakly mixing) and constructed a sensitive dynamical system whose $g$-fuzzification
is not sensitive for any $g\in D_{m}(I)$; giving a negative answer to Question \ref{Question 1}.
In this paper, we further investigate the relationships between some dynamical properties
(for example, transitivity, weakly mixing, mildly mixing, equicontinuity, uniform rigidity)
of $(K(X), T_{K})$ and $(\mathbb{F}_{0}(X), T_{F})$ through further developing the results in
\cite{Kupka2011}. In this study, we prove that  dynamical system is weakly mixing (resp.,
mildly mixing, weakly mixing and $\mathbf{a}$-transitive, equicontinuous, uniformly rigid) if and only if
the Zadeh's extension is transitive (resp., mildly mixing, $\mathbf{a}$-transitive, equicontinuous, uniformly rigid).

\medskip

This paper is organized as follows: in Section \ref{S-2}, some basic definitions and notations are introduced.
In Section \ref{S-3} and \ref{S-4}, the transitivity, the weakly mixing property and the mildly mixing
property of Zadeh's extension are studied. Then, in Sections \ref{S-5} and \ref{S-6}, some results on the
equicontinuity and the uniform rigidity are obtained.

\section{Basic definitions and notations}\label{S-2}

\subsection{Furstenberg family}
Let $\mathcal{P}$ be the
collection of all subsets of $\mathbb{Z}^{+}$. A collection
$\mathscr{F}\subset\mathcal{P}$ is called a {\it F\"{u}rstenberg family}
(briefly, a {\it family}) if it is hereditary upwards, i.e., $F_{1}\subset F_{2}$ and
$F_{1}\in\mathscr{F}$ imply $F_{2}\in\mathscr{F}$. A family $\mathscr{F}$
is {\it proper} if it is a proper subset of $\mathcal{P}$, i.e.,
neither empty nor the whole $\mathcal{P}$. It is easy to see
that $\mathscr{F}$ is proper if and only if $\mathbb{Z}^{+}\in \mathscr{F}$
and $\emptyset\notin \mathscr{F}$. Let $\mathscr{F}_{inf}$ be a
family of all infinite subsets of $\mathbb{Z}^+$.

\medskip

Given a family $\mathscr{F}$, define its {\it dual
family} as
\[
\kappa\mathscr{F}=\left\{F\in \mathcal{P}:
F\cap F'\neq \emptyset \text{ for all } F'\in \mathscr{F}\right\}=\left\{F\in \mathcal{P}:
\mathbb{Z}^+\setminus F\notin \mathscr{F}\right\}.
\]
Clearly, $\kappa\mathscr{F}_{inf}$ is the family of all cofinite subsets.
It is easy to check that $\kappa\mathscr{F}$ is a  family,
and is proper if and only if $\mathscr{F}$ is so. Given two  families $\mathscr{F}_1$
and $\mathscr{F}_2$, define $\mathscr{F}_{1}\cdot \mathscr{F}_{2}=\left\{F_{1}\cap F_{2}:
F_{1}\in \mathscr{F}_{1}, F_{2}\in \mathscr{F}_{2}\right\}$.
A  family $\mathscr{F}$ is {\it full} if $\mathscr{F}\cdot \kappa\mathscr{F}\subset
\mathscr{F}_{inf}$. A  family $\mathscr{F}$ is a {\it filter} if it is
proper and satisfies $\mathscr{F}\cdot\mathscr{F}\subset \mathscr{F}$; and it is
a {\it filterdual} if its dual family $\kappa\mathscr{F}$ is a filter.
For a  family $\mathscr{F}$, let $\mathscr{F}-\mathscr{F}=\{F-F: F\in \mathscr{F}\}$,
where $F-F=\{i-j: i, j\in F\}\cap \mathbb{Z^+}$.

\medskip

A subset $S$ of $\mathbb{Z}^+$ is {\it syndetic} if it has a bounded
gap, i.e., if there is $N \in \mathbb{N}$ such that $\{i, i+1, \ldots,
i+N\} \cap S \neq \emptyset$ for every $i \in \mathbb{Z}^+$; $S$ is
{\it thick} if it contains arbitrarily long runs of positive
integers, i.e., for every $n \in \mathbb{N}$ there exists some $a_n
\in \mathbb{Z}^+$ such that $\{a_n, a_n+1, \ldots, a_n+n\} \subset
S$. The set of all thick subsets of $\mathbb{Z}^{+}$ and all syndetic
subsets of $\mathbb{Z}^{+}$ are denoted by $\mathscr{F}_{t}$ and $\mathscr{F}_{s}$,
respectively. Clearly, they are both families.

\medskip

Let $\{p_{i}\}_{i=1}^{\infty}$ be an infinite sequence in $\mathbb{N}$ and
$$
FS(\{p_{i}\}_{i=1}^{\infty})=\left\{p_{i_1}+p_{i_2}+\cdots +p_{i_n}: 1\leq i_1<i_2<\cdots <i_n,\ n\in \mathbb{N}\right\}.
$$
A subset $A\subset \mathbb{Z}^+$ is an {\it IP-set} if it equals to some $FS(\{p_{i}\}_{i=1}^{\infty})$.
Denote the  family generated by all IP-sets by $\mathscr{F}_{IP}$. It follows from
Hindman's Theorem  \cite{Hin} that $\mathscr{F}_{IP}$ is a filterdual.

\medskip

For a  family $\mathscr{F}$, a dynamical system $(X, T)$
is called {\it $\mathscr{F}$-transitive} if $N(U, V)\in \mathscr{F}$
for every pair of nonempty open subsets $U, V\subset X$; and it is
{\it $\mathscr{F}$-mixing} if $(X \times X, T\times T)$ is $\mathscr{F}$-transitive.

\begin{lem}\label{F-mixing}\cite{A}
Let $(X, T)$ be a dynamical system and $\mathscr{F}$ be a full  family. Then,
the following statements are equivalent:
\begin{enumerate}[(1)]
\item $(X, T)$ is $\mathscr{F}$-mixing;
\item $\forall n\in \mathbb{N}$, $T^{(n)}$ is $\mathscr{F}$-mixing;
\item $(X, T)$ is weakly mixing
and $\mathscr{F}$-transitive.
\end{enumerate}
\end{lem}

\subsection{Topological dynamics}
For $U, V\subset X$, define the {\it return time set from $U$
to $V$} as $ N(U, V)=\{n\in \mathbb{Z}^{+}:
T^{n}(U)\cap V\neq \emptyset\}$. In particular,
$N(x, V)=\left\{n\in \mathbb{Z}^{+}:
T^{n}(x)\in V\right\}$ for $x\in X$.

\medskip

A dynamical system $(X, T)$ is
\begin{enumerate}[(1)]
\item {\it (topologically) weakly mixing}
if $(X \times X, T\times T)$ is transitive;
\item {\it topologically mixing} if for every pair of nonempty open subsets $U, V$ of $X$, $N(U, V)$ is cofinite,
i.e., there exists $m\in \mathbb{N}$ such that $[m, +\infty)\subset
N(U, V)$.
\end{enumerate}

It is well known that $(X, T)$ is transitive (resp., weakly mixing) if and only if
it is $\mathscr{F}_{inf}$-transitive (resp., $\mathscr{F}_{t}$-transitive) (see \cite{Fur-book}).
A point $x \in X$ is a {\it transitive point} of $T$
if its orbit $\mathrm{orb}(x, T):=\{x, T(x), T^{2}(x), \ldots \}$ is dense
in $X$. The set of all transitive points of $T$ is denoted by
$\mathrm{Tran}(T)$. It is well known that if $(X, T)$ is transitive,
then $\mathrm{Tran}(T)$ is a dense $G_{\delta}$-set.
The {\it $\omega$-limit set} of $x$ is the set of limit points of its orbit sequence
$\omega(x, T)=\bigcap_{m=0}^{+\infty}\overline{\{T^{n}(x): n\geq m\}}$. A
point $x\in X$ is a {\it recurrent point} of $T$ if $x\in \omega(x, T)$, i.e., there
exists $m_{i}\longrightarrow +\infty$ such that $T^{m_{i}}(x)\longrightarrow x$.
A well known result of Birkhoff states that every dynamical system admits a recurrent
point.

\begin{lem}\cite{Fur-book}\label{weakly-mixing}
Let $(X, T)$ be a dynamical system. Then, the following statements are equivalent:
\begin{enumerate}[(1)]
\item $(X, T)$ is weakly mixing;
\item For any pair of nonempty open subsets $U, V$ of $X$, $N(U, U)\cap N(U, V)\neq \emptyset$;
\item For any pair of nonempty open subsets $U, V$ of $X$, $N(U, V)\cap N(V, V)\neq \emptyset$;
\item $(X, T)$ is $\mathscr{F}_{t}$-transitive.
\end{enumerate}
\end{lem}

Recently, Moothathu~\cite{Moothathu-2010} introduced the notion of multi-transitivity.
A dynamical system $(X, T)$ is called \emph{multi-transitive} if for any $n\in\mathbb{N}$,
the product system $(X^n, T\times T^2\times\dotsb\times T^n)$ is transitive.
He also proved that a minimal system is multi-transitive if and only if it is weakly mixing and
asked whether there are implications between the multi-transitivity and the weak mixing property
for general (not necessarily minimal) systems. Then, Kwietniak and Oprocha~\cite{D-Kwietniak-P-Oprocha-2010}
showed that in general there is no connection between the multi-transitivity and the weakly mixing property
by constructing examples of weakly mixing but non-multi-transitive
and multi-transitive but non-weakly mixing systems.
To generalize the concept of multi-transitivity, Chen et al. \cite{CLL14} introduced
the notion of multi-transitivity with respect to a vector.
Let $\mathbf{a}=(a_{1}, a_{2}, \dotsc, a_n)$ be a vector in $\mathbb{N}^n$.
A dynamical system $(X, T)$ is \emph{multi-transitive with respect to the vector $\mathbf{a}$}
(briefly, \emph{$\mathbf{a}$-transitive}) if the product system $(X^{n},  T^{(\mathbf{a})})$ is transitive,
where $T^{(\mathbf{a})}=T^{a_{1}}\times T^{a_{2}}\times\cdots \times T^{a_{r}}$.

\begin{lem}\label{a-tran}\cite{Wu}
Let $(X, T)$ be a dynamical system and $\mathbf{a}=(a_{1}, \ldots, a_{n})\in \mathbb{N}^{n}$. Then,
the following statements are equivalent:
\begin{enumerate}[(1)]
  \item\label{4-2-1} $(X, T)$ is weakly mixing and $\mathbf{a}$-transitive;
  \item\label{4-2-2} For any $m\in \mathbb{N}$, $(X^{m}, T^{b_{1}}\times T^{b_{2}}\times \cdots \times T^{b_{m}})$
  is transitive, where $b_{i}\in \{a_{i}: 1\leq i\leq n\}$, $i=1, 2, \ldots, m$.
\end{enumerate}
\end{lem}

A dynamical system $(X, T)$ is {\it weakly disjoint} with another dynamical system $(Y, S)$ if their
product system $(X\times Y, T\times S)$ is transitive. A dynamical system is {\it mildly mixing} if
it is weakly disjoint with every transitive system. Huang and Ye \cite{HY} proved that a dynamical system
is mildly mixing if and only if it is $\kappa(\mathscr{F}_{IP}-\mathscr{F}_{IP})$-transitive.

\medskip

A dynamical system $(X, T)$ is {\it equicontinuous} if for any $\varepsilon>0$, there exists $\delta>0$ such that for
any $x, y\in X$ with $d(x, y)< \delta$ and any $n\in \mathbb{Z}^{+}$,
$d(T^{n}(x), T^{n}(y))<\varepsilon$.

\subsection{Hyperspace $K(X)$}
Let $K(X)$ be the hyperspace on $X$, i.e., the space of all nonempty
closed subsets of $X$ with the Hausdorff metric $d_{H}$ defined by
\[
d_{H}(A, B)=\max\left\{\max_{x\in A}\min_{y\in B}d(x, y), \max_{y\in B}\min_{x\in A}d(x, y)\right\},
\quad \forall A, B\in K(X).
\]
It is known that
$(K(X), d_{H})$ is also a compact metric space (see \cite{Illanes1999}).
The system $(X, T)$ induces naturally
a set-valued dynamical system $(K(X), T_{K})$,
where $T_{K}: K(X)\longrightarrow K(X)$ is defined as
$T_{K}(A)=T(A)$ for any $A\in K(X)$.
For any finite collection $A_{1}, \ldots, A_{n}$ of nonempty subsets of $X$, take
\[
\langle A_{1}, \ldots, A_{n}\rangle=\left\{A\in K(X): A\subset \bigcup_{i=1}^{n}A_{i},
A\cap A_{i}\neq \emptyset \text{ for all } i=1, \ldots, n\right\}.
\]
It follows from \cite{Illanes1999} that the topology on $K(X)$
given by the metric $d_{H}$ is same as the Vietoris or finite
topology, which is generated by a basis consisting of all sets of the following form,
\[
\langle U_{1}, \ldots, U_{n}\rangle,
\text{ where } U_{1}, \ldots, U_{n} \text{ are an arbitrary finite collection of nonempty open subsets of
 } X.
\]
Under this topology
$\mathcal{F}(X)$, the set of all finite subsets of $X$,
is dense in $K(X)$.

\subsection{Zadeh's extension}

Let $I= [0, 1]$. A {\it fuzzy set $A$} in space $X$ is a function $A: X\longrightarrow I$.
Given a fuzzy set $A$,
its {\it $\alpha$-cuts} (or {\it $\alpha$-level sets}) $[A]_{\alpha}$
and {\it support} $\mathrm{supp}(A)$ are defined respectively by
\[
[A]_{\alpha}=\{x\in X: A(x)\geq \alpha\}, \quad \forall \alpha\in I,
\]
and
\[
\mathrm{supp}(A)=\overline{\left\{x\in X: A(x)>0\right\}}.
\]

Let $\mathbb{F}(X)$ denote the set of all upper semicontinuous fuzzy sets defined on $X$
and set
$$
\mathbb{F}^{\geq \lambda}(X)=\left\{A\in \mathbb{F}(X): A(x)\geq \lambda \text{ for some }
x\in X\right\},
$$
$$
\mathbb{F}^{=\lambda}(X)=\left\{A\in \mathbb{F}(X): \max\{A(x): x\in X\}=\lambda\right\}.
$$
Especially, let $\mathbb{F}^{=1}(X)$ denote the system of all normal fuzzy sets on $X$.

Define $\emptyset_{X}$ as the empty fuzzy set ($\emptyset_{X}\equiv 0$) in $X$,
and $\mathbb{F}_{0}(X)$ as the set of all nonempty upper semicontinuous fuzzy sets.
Since the Hausdorff metric $d_{H}$ is measured only between two nonempty closed subsets
in $X$, one can consider the following extension of the Hausdorff metric:
$$
d_{H}(\emptyset, \emptyset)=0 \text{ and } d_{H}(\emptyset, A)=d_{H}(A, \emptyset)=
\mathrm{diam}(X), \quad \forall A\in {K}(X).
$$
Under this Hausdorff metric, one can define a {\it levelwise metric} $d_{\infty}$ on $\mathbb{F}(X)$ by
$$
d_{\infty}(A, B)=
\sup\left\{d_{H}([A]_{\alpha}, [B]_{\alpha}): \alpha\in (0, 1]\right\}, \quad \forall A, B\in \mathbb{F}(X).
$$
It is well known that the spaces $(\mathbb{F}(X), d_{\infty})$
and $(\mathbb{F}^1(X), d_{\infty})$ are complete, but not compact and not
separable (see \cite{Kupka2011-1} and references therein).

A fuzzy set $A\in \mathbb{F}(X)$ is {\it piecewise constant} if there exists a finite
number of sets $D_{i}\subset X$ such that $\bigcup \overline{D_{i}}=X$ and $A|_{\mathrm{int}\overline{D_{i}}}$
is constant. In this case, a piecewise constant $A$ can be represented
by a strictly decreasing sequence of closed subsets $\{A_{1}, A_{2}, \ldots, A_{k}\}\subset K(X)$
and a strictly increasing sequence of reals $\{\alpha_{1}, \alpha_{2}, \ldots, \alpha_{k}=\max\{A(x): x\in X\}\}\subset (0, 1]$ if
\[
[A]_{\alpha}=A_{i+1}, \text{ whenever } \alpha\in (\alpha_{i}, \alpha_{i+1}].
\]

\begin{rmk}\label{remark-piecewise}
Fix any two piecewise constants $A, B\in \mathbb{F}(X)$ which are represented by
strictly decreasing sequences of closed subsets $\{A_{1}, A_{2}, \ldots, A_{k}\}$,
$\{B_{1}, B_{2}, \ldots, B_{s}\}\subset K(X)$ and strictly increasing sequences of
reals $\{\alpha_{1}, \alpha_{2}, \ldots, \alpha_{k}\}$, $\{\beta_{1}, \beta_{2}, \ldots, \beta_{s}\}\subset (0, 1]$
with
$$
[A]_{\alpha}=A_{i+1}, \quad \forall \alpha\in (\alpha_{i}, \alpha_{i+1}] \text{ and }
[B]_{\alpha}=B_{i+1}, \quad \forall \beta\in (\beta_{i}, \beta_{i+1}],
$$
respectively. Arrange all reals $\alpha_{1}, \alpha_{2}, \ldots, \alpha_{k}, \beta_{1}, \beta_{2}, \ldots, \beta_{s}$
by the natural order `$<$' and denote them
by $\gamma_{1}, \gamma_{2}, \ldots, \gamma_{n}$ $(n\leq k+s)$. Then, it can be verified
that for any $1\leq t<n$, there exist $1\leq i\leq k$ and $1\leq j\leq s$ such that for
any $\gamma\in (\gamma_{t}, \gamma_{t+1}]$,
$$
[A]_{\gamma}=A_{i} \text{ and } [B]_{\gamma}=B_{j}.
$$
This implies that there exist (not necessarily strictly) decreasing sequences of closed subsets
$\{C_{1}, C_{2}, \ldots, C_{n}\}$, $\{D_{1}, D_{2}, \ldots, D_{n}\}\subset K(X)$
and a strictly increasing sequence of reals $\gamma_{1}, \gamma_{2}, \ldots, \gamma_{n}
\subset (0, 1]$ such that
$$
[A]_{\gamma}=C_{i+1} \text{ and } B_{\gamma}=D_{i+1}, \text{ whenever }
\gamma\in (\gamma_{i}, \gamma_{i+1}].
$$
\end{rmk}

To generalize the concept of Zadeh's extension,
Kupka \cite{Kupka2011-1} introduced the notion of $g$-fuzzification.

\medskip

{\it Zadeh's extension} (also called {\it usual fuzzification}) of a dynamical system
$(X, T)$ is a map $T_{F}: \mathbb{F}(X)\longrightarrow \mathbb{F}(X)$
defined by
\[
T_{F}(A)(x)=\sup\left\{A(y): y\in T^{-1}(x)\right\},
\quad \forall A\in \mathbb{F}(X), \forall
x\in X.
\]

\begin{lem}\label{piecewise}\cite[Lemma 1, Remark 1, Theorem 1]{Kupka2011}
Let $(X, T)$ be a dynamical system and $\lambda\in (0, 1]$. Then,
\begin{enumerate}[(1)]
\item the set of piecewise constants is dense in $\mathbb{F}(X)$,
$\mathbb{F}^{\geq \lambda}(X)$ and $\mathbb{F}^{=\lambda}(X)$.
\item $(\mathbb{F}^{=\lambda}(X), T_{F}|_{\mathbb{F}^{=\lambda}(X)})$
is periodically dense in $\mathbb{F}^{=\lambda}(X)$ if and only if $(K(X), T_{K})$
is periodically dense in $K(X)$.
\end{enumerate}
\end{lem}

\begin{rmk}\label{Remark 3}
For any $A, B\in K(X)$, any $n\in \mathbb{N}$ and any $\lambda\in (0, 1]$,
\begin{equation}\label{e-4}
d_{H}(T_{K}^{n}(A), T_{K}^{n}(B))=d_{\infty}(T_{F}^{n}(\lambda\cdot \chi_{A}),
T_{F}^{n}(\lambda\cdot \chi_{B})).
\end{equation}
This shows that for any fixed $\lambda\in (0, 1]$, the subsystem $(\mathbb{F}_{\lambda\cdot \chi}:=
\{\lambda \cdot \chi_{A}\in \mathbb{F}_{0}(X):
A\in K(X)\}, T_{F}|_{\mathbb{F}_{\lambda\cdot \chi}})$ is topologically conjugated to
$(K(X), T_{K})$.
\end{rmk}

Let $D_{m}(I)$ be the set of all nondecreasing right-continuous functions
$g: I\longrightarrow I$ with $g(0)=0$ and $g(1)=1$. For a dynamical system $(X, T)$
and for any $g\in D_{m}(I)$, define a map $T_{F}^{g}: \mathbb{F}(X)
\longrightarrow \mathbb{F}(X)$ by
\[
T_{F}^{g}(A)(x)=\sup\left\{g(A(y)): y\in T^{-1}(x)\right\},
\quad \forall A\in \mathbb{F}(X), \forall
x\in X,
\]
which is called the {\it $g$-fuzzification} of the dynamical system $(X, T)$.
Clearly, $T_{F}=T_{F}^{\mathrm{id}_{I}}$, where $\mathrm{id}_{I}$ is the
identity map defined on $I$.

Also, define the {\it $\alpha$-cut} $[A]_{\alpha}^{g}$ of a fuzzy set
$A\in \mathbb{F}(X)$ with respect to $g\in D_{m}(I)$ by
\[
[A]_{\alpha}^{g}=\left\{x\in \mathrm{supp}(A): g(A(x))\geq \alpha\right\}.
\]

For any $g\in D_{m}(I)$, the right-continuity of $g$
implies that $\min g^{-1}([x, 1])$ exists for any $x\in [0, 1]$.
Since $g$ is nondecreasing, $\min g^{-1}([x, 1])>0$
holds for any $x\in (0, 1]$. Define $\xi_{g}: [0, 1]\longrightarrow [0,1]$
by $\xi_{g}(x)=\min g^{-1}([x, 1])$ for any $x\in [0, 1]$. Clearly, $\xi_{g}$
is nondecreasing. Recently, we \cite{WC-2} proved the following result:

\begin{lem}\cite{WC-2}\label{Lemma}
Let $(X, T)$ be a dynamical system, $g\in D_{m}(I)$ and $T_{F}^{g}$
be the $g$-fuzzification of $T$. Then, for any $n\in \mathbb{N}$, any $A\in \mathbb{F}(X)$
and any $\alpha\in (0, 1]$, $\left[(T_{F}^{g})^{n}(A)\right]_{\alpha}
=T^{n}([A]_{\xi_{g}^{n}(\alpha)})$. In particular, $\left[T_{F}^{n}(A)\right]_{\alpha}
=T^{n}([A]_{\alpha})$.
\end{lem}

\section{Transitivity of $(\mathbb{F}^{=\lambda}(X), T_{F})$}\label{S-3}

Banks \cite{Banks2005} proved the following result on the transitivity
of $(K(X), T_{K})$.

\begin{lem}\cite[Theorem 2]{Banks2005}\label{K(X)}
Let $(X, T)$ be a dynamical system. Then, the following statements are equivalent:
\begin{enumerate}[(1)]
\item $(X, T)$ is weakly mixing;
\item $(K(X), T_{K})$ is weakly mixing;
\item $(K(X, T_{K})$ is transitive.
\end{enumerate}
\end{lem}

Inspired by Lemma \ref{K(X)}, this section is devoted to studying the transitivity of Zadeh's extension $(\mathbb{F}^{=\lambda}(X), T_{F}|_{\mathbb{F}^{=\lambda}(X)})$.
In particular, it is proved that both the transitivity and the weakly mixing property of
$(\mathbb{F}^{=\lambda}(X), T_{F}|_{\mathbb{F}^{=\lambda}(X)})$
are equivalent to the weakly mixing property of $(X, T)$ (see Theorem \ref{tran-F(X)}).

\begin{thm}
Let $(X, T)$ be a dynamical system and $C$
be an invariant closed subset of $\mathbb{F}_{0}(X)$. If $(C, T_{F}|_{C})$
is transitive, then for any $A, B\in C$, $\max\{A(x): x\in X\}=\max\{B(x): x\in X\}$,
i.e., there exists $\lambda\in [0, 1]$ such that $C\subset \mathbb{F}^{=\lambda}(X)$.
\end{thm}

\pf
Suppose that there exist $A, B\in C$ such that
$$
\xi:=\max\{A(x): x\in X\}<\max\{B(x): x\in X\}:=\eta.
$$
Choose $a=\min\left\{\frac{\eta-\xi}{4}, \frac{\mathrm{diam}(X)}{4}\right\}>0$ and set
$U=B_{d_{\infty}}(A, a)\cap C$ and $V=B_{d_{\infty}}(B, a)\cap C$. Clearly, $U$ and $V$
are nonempty open subsets of $C$. For any $D\in U$,
$$
a>d_{\infty}(A,D)\geq d_{H}([A]_{a}, [D]_{a})=d_{H}(\emptyset, [D]_{a}),
$$
implying that $[D]_{a}=\emptyset$. Then, for any $n\in \mathbb{Z}^{+}$,
\begin{eqnarray*}
d_{H}([T_{F}^{n}(D)]_{a}, [B]_{a})&=&d_{H}(T_{K}^{n}([D]_{a}), [B]_{a})\\
&=&d_{H}(\emptyset, [B]_{a}) \quad (\text{as } [B]_{a}\neq \emptyset)\\
&=& \mathrm{diam}(X).
\end{eqnarray*}
This implies that $T_{F}^{n}(U)\cap V=\emptyset$, i.e., $(C, T_{F}|_{C})$
is not transitive, which is a contradiction.
\epf



\begin{lem}\label{lemma 1}
Let $(X, T)$ be a dynamical system and $\lambda\in (0, 1]$. If
$(\mathbb{F}^{=\lambda}(X), T_{F}|_{\mathbb{F}^{=\lambda}(X)})$ is transitive, then
$(K(X), T_{K})$ is weakly mixing.
\end{lem}

\pf Applying Lemma \ref{K(X)}, it suffices to prove that $T_{K}$
is transitive.

For any pair of nonempty open subsets $U, V$ of $K(X)$,
there exist $A\in U$, $B\in V$ and $0<\delta<\frac{\mathrm{diam}(X)}{2}$ such that
$B_{d_{H}}(A, \delta):=\{C\in K(X): d_{H}(C, A)<\delta\}\subset U$ and
$B_{d_{H}}(B, \delta)\subset V$. Noting that $U_{1}:=B_{d_{\infty}}(\lambda\cdot\chi_{A}, \delta)\cap
\mathbb{F}^{=\lambda}(X)$ and $V_{1}:=B_{d_{\infty}}(\lambda\cdot \chi_{B}, \delta)\cap \mathbb{F}^{=\lambda}(X)$
are nonempty open subsets of $\mathbb{F}^{=\lambda}(X)$,
since $T_{F}|_{\mathbb{F}^{=\lambda}(X)}$ is transitive, there exists $n\in \mathbb{Z}^{+}$ such that
\[
T_{F}^{n}(U_{1})\cap V_{1}\neq \emptyset.
\]
Then, there exists a point $E\in U_1$ such that
$T_{F}^{n}(E)\in V_{1}$. This implies that
\begin{equation}\label{1}
d_{H}\left(\left[T_{F}^{n}(E)\right]_{\lambda}, [\lambda\cdot \chi_B]_{\lambda}\right)=d_{H}(T_{K}^{n}([E]_{\lambda}), B)<\delta.
\end{equation}
Since $E\in B_{d_{\infty}}(\lambda\cdot \chi_A, \delta)$,
it can be verified that
$$
d_{H}(A, [E]_{\lambda})<\delta.
$$
Clearly,
$$
[E]_{\lambda}\neq \emptyset.
$$
Then,
$$
[E]_{\lambda}\in B_{d_{H}}
(A, \delta)\subset U.
$$
Combining this with \eqref{1},
it follows that
\[
T_{K}^{n}([F]_{\lambda})\in T_{K}^{n}(U)\cap V\neq \emptyset.
\]
\epf

\begin{lem}\label{lemma 2}
Let $(X, T)$ be a dynamical system and $\lambda\in (0, 1]$.
If $(K(X), T_{K})$ is weakly mixing, then
$(\mathbb{F}^{=\lambda}(X), T_{F}|_{\mathbb{F}^{=\lambda}(X)})$ is weakly mixing.
\end{lem}

\pf
Given any pair of nonempty open subsets $U, V$ of $\mathbb{F}^{=\lambda}(X)$,
applying Lemma \ref{piecewise} implies that there exist piecewise constants
$P\in U$, $Q\in V$ and $\delta>0$ such that $B_{d_{H}}(P, \delta)\cap \mathbb{F}^{=\lambda}(X)
\subset U$ and $B_{d_{H}}(Q, \delta)\cap \mathbb{F}^{=\lambda}(X)\subset V$. Since $P$ and $Q$
are piecewise constants and $P, Q\in \mathbb{F}^{=\lambda}(X)$, it follows from
Remark \ref{remark-piecewise} that there exist
strictly increasing sequence of reals
$\{\alpha_{1}, \ldots, \alpha_{k}=\lambda\}\subset [0, 1]$
and decreasing sequences of closed subsets $\{C_{1}, \ldots, C_{k}\}\subset K(X)$,
$\{D_{1}, \ldots, D_{k}\}\subset K(X)$ such that
\[
[P]_{\alpha}=C_{i+1}, \ [Q]_{\alpha}=D_{i+1}, \text{ whenever } \alpha\in (\alpha_{i}, \alpha_{i+1}].
\]
Noting that $B_{d_{H}}(C_{1}, \frac{\delta}{2}), \ldots, B_{d_{H}}(C_{k}, \frac{\delta}{2})$ and
$B_{d_{H}}(C_{1}, \frac{\delta}{2}), \ldots, B_{d_{H}}(C_{k}, \frac{\delta}{2})$ are nonempty
open subsets of $K(X)$, since $T_{K}$ is weakly mixing, it
follows that there exists $n\in \mathbb{Z}^{+}$ such that for any $1\leq i\leq k$,
\[
T_{K}^{n}\left(B_{d_{H}}
\left(C_{i}, \frac{\delta}{2}\right)\right)\cap B_{d_{H}}
\left(C_{i}, \frac{\delta}{2}\right) \neq \emptyset,
\]
and
\[
T_{K}^{n}\left(B_{d_{H}}
\left(C_{i}, \frac{\delta}{2}\right)\right)\cap B_{d_{H}}
\left(D_{i}, \frac{\delta}{2}\right) \neq \emptyset,
\]
implying that there exist $G_{i}, E_{i}\in B_{d_{H}}(C_{i}, \frac{\delta}{2})$ such that
$$
T_K^{n}(G_{i})\in B_{d_{H}}\left(C_{i}, \frac{\delta}{2}\right) \text{ and }
T_{K}^{n}(E_{i})\in B_{d_{H}}\left(D_{i}, \frac{\delta}{2}\right).
$$
Define respectively
two fuzzy sets $G: X\longrightarrow I$ and $E: X\longrightarrow I$ by
\[
[G]_{\alpha}=\bigcup_{j=i+1}^{k}G_{j} \text{ and }
[E]_{\alpha}=\bigcup_{j=i+1}^{k}E_{j}, \quad \forall \alpha\in (\alpha_{i}, \alpha_{i+1}].
\]
Clearly, $G, E\in \mathbb{F}^{=\lambda}(X)$. Since $G_{i}, E_{i}\in B_{d_{H}}(C_{i}, \frac{\delta}{2})$, it can be verified
that for any $\alpha\in (\alpha_{i}, \alpha_{i+1}]$,
\[
d_{H}([G]_{\alpha}, [P]_{\alpha})=d_{H}\left(\bigcup_{j=i+1}^{k}G_{j}, C_{i+1}\right)<\frac{\delta}{2},
\]
and
\[
d_{H}([E]_{\alpha}, [P]_{\alpha})<d_{H}\left(\bigcup_{j=i+1}^{k}E_{j}, C_{i+1}\right)<\frac{\delta}{2}.
\]
This implies that
$$
G, E\in B_{d_{\infty}}(P, \delta)\cap \mathbb{F}^{=\lambda}(X)\subset U.
$$
Meanwhile,
since $T_{K}^{n}(G_{i})\in B_{d_{H}}(C_{i}, \frac{\delta}{2})$ and
$T_{K}^{n}(E_{i})\in B_{d_{H}}(D_{i}, \frac{\delta}{2})$, it can be verified
that for any $\alpha\in (\alpha_{i}, \alpha_{i+1}]$,
\[
d_{H}([T_{F}^{n}(G)]_{\alpha}, [P]_{\alpha})=
d_{H}(T_{K}^{n}([G]_{\alpha}), [P]_{\alpha})=d_{H}\left(\bigcup_{j=i+1}^{k}T_{K}^{n}(G_{j}), C_{i+1}\right)<\frac{\delta}{2},
\]
and
\[
d_{H}([T_{F}^{n}(E)]_{\alpha}, [Q]_{\alpha})=
d_{H}(T_{K}^{n}([E]_{\alpha}), [Q]_{\alpha})=d_{H}\left(\bigcup_{j=i+1}^{k}T_{K}^{n}(E_{j}), D_{i+1}\right)<\frac{\delta}{2}.
\]
Then,
$$
\left(T_{F}|_{\mathbb{F}^{=\lambda}(X)}\right)^{n}(G)\in B_{d_{H}}(P, \delta)\cap \mathbb{F}^{=\lambda}(X)\subset U
$$
and
$$
\left(T_{F}|_{\mathbb{F}^{=\lambda}(X)}\right)^{n}(E)\in B_{d_{H}}(Q, \delta)\cap \mathbb{F}^{=\lambda}(X)\subset V.
$$
Therefore,
\[
\left(T_{F}|_{\mathbb{F}^{=\lambda}(X)}\right)^{n}(G)\in T_{F}^{n}(U)\cap U\neq \emptyset \text{ and }
\left(T_{F}|_{\mathbb{F}^{=\lambda}(X)}\right)^{n}(E)\in T_{F}^{n}(U)\cap V\neq \emptyset.
\]
This, together with Lemma \ref{weakly-mixing}, implies that
$(\mathbb{F}^{=\lambda}(X), T_{F})$ is weakly mixing.
\epf

\begin{thm}\label{tran-F(X)}
Let $(X, T)$ be a dynamical system. Then, the following statements are equivalent:
\begin{enumerate}[(1)]
\item\label{1-1} $(X, T)$ is weakly mixing;
\item\label{1-2} $(K(X), T_{K})$ is transitive;
\item\label{1-3} $(K(X), T_{K})$ is weakly mixing;
\item\label{1-4} $\forall \lambda\in (0, 1]$, $(\mathbb{F}^{=\lambda}(X), T_{F}|_{\mathbb{F}^{=\lambda}(X)})$ is transitive;
\item\label{1-5} $\forall \lambda\in (0, 1]$, $(\mathbb{F}^{=\lambda}(X), T_{F}|_{\mathbb{F}^{=\lambda}(X)})$ is weakly mixing.
\end{enumerate}
\end{thm}
\pf
Applying Lemma \ref{K(X)} implies that (\ref{1-1})$\Longleftrightarrow$
(\ref{1-2})$\Longleftrightarrow$(\ref{1-3}). It follows from Lemma \ref{lemma 1}
and Lemma \ref{lemma 2} that (\ref{1-4})$\Longrightarrow$(\ref{1-2})
$\Longrightarrow$(\ref{1-5})$\Longrightarrow$(\ref{1-4}).
\epf


\begin{coy}\label{Corollary 3.1}
Let $(X, T)$ be a dynamical system. Then, the following statements are equivalent:
\begin{enumerate}[(1)]
\item $({K}(X), T_{K})$ is Devaney chaotic;
\item $\forall \lambda\in (0, 1]$, $(\mathbb{F}^{=\lambda}(X), T_{F}|_{\mathbb{F}^{=\lambda}(X)})$ is Devaney chaotic;
\item $\forall \lambda\in (0, 1]$, $(\mathbb{F}^{=\lambda}(X), T_{F}|_{\mathbb{F}^{=\lambda}(X)})$ is Devaney chaotic.
\end{enumerate}
\end{coy}
\pf
It follows immediately from Theorem \ref{tran-F(X)} and Lemma \ref{piecewise}.
\epf

\begin{rmk}
\begin{enumerate}[(1)]
\item In \cite[Proposition 3.4, Proposition 4.2]{Lan2012}, Lan et al. proved that
$(X, T)$ is weakly mixing if and only if $(\mathbb{F}^{=1}(X), T_{F}|_{\mathbb{F}^{=1}(X)})$ is weakly mixing. However,
their proof is not correct. Because the proof of \cite[Proposition 4.2]{Lan2012}
is based on \cite[Proposition 4.1 (4)]{Lan2012} which claimed that for any nonempty
open subset $U\subset \mathbb{F}^{=1}(X)$,
\[
r(U):=\left\{A\in K(X):
\exists u\in U \text{ such that } A\subset [u]_{0}\right\}
\]
is also a nonempty open subset of $X$. It is clear that $r(U)$ is not a open subset
of $X$, because $K(X)\supset r(U)\nsubseteq X$. Meanwhile, noting that for any $u\in U$,
$[u]_{0}=X$, it is easy to see that $r(U)=K(X)$.

\item Kupka \cite{Kupka2011} obtained that the Devaney's chaoticity of
$(X, T)$ does not imply the same of $(\mathbb{F}^{=1}(X), T_{F}|_{\mathbb{F}^{=1}(X)})$.
According to \cite[Remark 2.4]{Li2014}, there exists a dynamical system $(X, T)$ such that
$(K(X), T_{K})$ is Devaney chaotic, while $(X, T)$ is not Devaney chaotic,
showing that the answer to \cite[Q3]{Liao2006} is negative.
This, together with Corollary \ref{Corollary 3.1}, shows that
the Devaney's chaoticity of $(\mathbb{F}^{=1}(X), T_{F}|_{\mathbb{F}^{=1}(X)})$
does not imply the Devaney's chaoticity of $(X, T)$.

\item Applying Theorem \ref{tran-F(X)} and Corollary \ref{Corollary 3.1} yields
that \cite[Theorem 2, Proposition 1, Theorem 3, Proposition 2, Theorem 4]{Kupka2011} holds
trivially and that the converses of \cite[Proposition 2, Theorem 4]{Kupka2011} are true.
\end{enumerate}
\end{rmk}

About the weakly mixing property of dynamical systems, Liao et al. \cite{Liao2006} provided the following question:

\begin{qun}\label{Q-1}\cite{Liao2006}
Which systems, besides $T_{K}$, have the equivalence between the transitivity
and the weakly mixing property?
\end{qun}

As a partial answer to Question \ref{Q-1}, applying Theorem \ref{tran-F(X)},
we know that the Zadeh's extension restricted on the space of normal fuzzy
sets has the equivalence between the transitivity and the weakly mixing property.

\begin{lem}\label{F-mixing-K(X)}
Let $(X, T)$ be a dynamical system and $\mathscr{F}$ be a full family. Then,
the following statements are equivalent:
\begin{enumerate}[(1)]
\item $(X, T)$ is $\mathscr{F}$-mixing;
\item $(K(X), T_{K})$ is $\mathscr{F}$-transitive;
\item $(K(X), T_{K})$ is $\mathscr{F}$-mixing.
\end{enumerate}
\end{lem}

Slightly modifying the proofs of Lemma \ref{lemma 1} and Lemma \ref{lemma 2},
applying Lemma \ref{F-mixing} and Lemma \ref{F-mixing-K(X)}, it is not difficult to prove the following.

\begin{thm}\label{mixing}
Let $(X, T)$ be a dynamical system and $\lambda\in (0, 1]$.
Then, $(X, T)$ is mixing if and only if $(K(X), T_{K})$ is
mixing if and only if $(\mathbb{F}^{=\lambda}(X), T_{F}|_{\mathbb{F}^{=\lambda}(X)})$ is mixing.
\end{thm}


\begin{coy}
Let $(X, T)$ be a dynamical system and $\mathscr{F}$ be a full family.
Then, the following statements are equivalent:
\begin{enumerate}[(1)]
\item $(X, T)$ is $\mathscr{F}$-mixing;
\item $(K(X), T_{K})$ is $\mathscr{F}$-transitive;
\item $(K(X), T_{K})$ is $\mathscr{F}$-mixing;
\item $\forall \lambda\in (0, 1]$,
$(\mathbb{F}^{=\lambda}(X), T_{F}|_{\mathbb{F}^{=\lambda}(X)})$ is $\mathscr{F}$-transitive;
\item $\forall \lambda\in (0, 1]$,
$(\mathbb{F}^{=\lambda}(X), T_{F}|_{\mathbb{F}^{=\lambda}(X)})$ is $\mathscr{F}$-mixing;
\item $\exists \lambda\in (0, 1]$,
$(\mathbb{F}^{=\lambda}(X), T_{F}|_{\mathbb{F}^{=\lambda}(X)})$ is $\mathscr{F}$-transitive;
\item $\exists \lambda\in (0, 1]$,
$(\mathbb{F}^{=\lambda}(X), T_{F}|_{\mathbb{F}^{=\lambda}(X)})$ is $\mathscr{F}$-mixing.
\end{enumerate}
\end{coy}

\section{Mildly mixing property and $\mathbf{a}$-transitivity of $(\mathbb{F}^{=\lambda}(X), T_{F})$}\label{S-4}

Bauer and Sigmund \cite{BS75} proved the equivalence of the mildly mixing property between
$(X, T)$ and $(K(X), T_{K})$.
\begin{lem}\cite[Theorem 1, Proposition 2]{BS75}\label{T=T_K}
A dynamical system $(X, T)$ is mildly mixing if and only if $(K(X), T_{K})$ is mildly mixing.
\end{lem}

Similarly to the proof of \cite[Proposition 7]{Kupka2014}, it can be verified that the
following result holds.

\begin{lem}\label{T_F-->T}
Let $(X, T)$ be a dynamical system and $\lambda\in (0, 1]$. If $(\mathbb{F}^{=\lambda}(X), T_{F}|_{\mathbb{F}^{=\lambda}(X)})$
is mildly mixing, then $(X, T)$ is mildly mixing.
\end{lem}

\begin{lem}\label{T_K-->T_F}
If $(K(X), T_{K})$ is mildly mixing, then for any $\lambda\in (0, 1]$,
$(\mathbb{F}^{=\lambda}(X), T_{F}|_{\mathbb{F}^{=\lambda}(X)})$ is mildly mixing.
\end{lem}

\pf
It suffices to check that for any transitive system $(Y, S)$, $(\mathbb{F}^{=\lambda}(X)
\times Y, T_{F}|_{\mathbb{F}^{=\lambda}(X)}\times S)$ is transitive.

For any pair of nonempty open subsets $W, V$ of $\mathbb{F}^{=\lambda}(X)
\times Y$, it follows from Lemma \ref{piecewise} and Remark \ref{remark-piecewise}
that there exist piecewise constants $A, B\in \mathbb{F}^{=\lambda}(X)$ which are
represented by decreasing sequences of closed subsets $\{A_{1}, A_{2}, \ldots, A_{k}\}$,
$\{B_{1}, B_{2}, \ldots, B_{k}\}\subset K(X)$ and a strictly increasing sequence
of reals $\{\alpha_{1}, \alpha_{2}, \ldots, \alpha_{k}=\lambda\}\subset (0, 1]$ such that
\[
[A]_{\alpha}=A_{i+1}, \ [B]_{\alpha}=B_{i+1}, \text{ whenever } \alpha\in (\alpha_{i}, \alpha_{i+1}],
\]
$y_{1}, y_{2}\in Y$ and $\delta>0$ such that
$$
\left[B_{d_{\infty}}(A, \delta)\cap\mathbb{F}^{=\lambda}(X)\right]\cap B(y_{1}, \delta)
\subset W,
$$
and
$$
\left[B_{d_{\infty}}(B, \delta)\cap\mathbb{F}^{=\lambda}(X)\right]\cap B(y_{2}, \delta)
\subset V.
$$
Clearly, $(\underbrace{K(X)\times \cdots \times K(X)}\limits_{k}\times Y,
\underbrace{T_{K}\times \cdots \times T_{K}}\limits_{k}\times S)$
is transitive, as $T_{K}$ is mildly mixing. This implies that there exists
$n\in \mathbb{Z}^+$ such that for any $1\leq i\leq k$,
$$
T_{K}^{n}(B_{d_{H}}(A_{i}, \delta))\cap B_{d_{H}}(B_{i}, \delta)\neq \emptyset,
$$
and
$$
S^{n}(B(y_1, \delta))\cap B(y_2, \delta)\neq \emptyset.
$$
Then there exist $C_{i}\in B_{d_{H}}(A_{i}, \delta)$ and $y\in B(y_{1}, \delta)$ such that
$$
d_{H}(T_{K}^{n}(C_{i}), B_{i})<\delta \text{ and } d(S^{n}(y), y_{2})<\delta.
$$
Take a piecewise constant $C\in \mathbb{F}^{=\lambda}(X)$ as
$$
[C]_{\alpha}=\bigcup_{j=i+1}^{k}C_{j}, \text{ whenever } \alpha\in (\alpha_{i}, \alpha_{i+1}].
$$
It can be verified that the following statements hold:
\begin{enumerate}[(i)]
\item $d_{H}\left(\bigcup_{j=i+1}^{k}C_{j}, A_{i+1}\right)=d_{H}\left(\bigcup_{j=i+1}^{k}C_{j}, \bigcup_{j=i+1}^{k}A_{j}\right)<\delta$;
\item $d_{H}\left(T_{K}^{n}(\bigcup_{j=i+1}^{k}C_{j}), B_{i+1}\right)=d_{H}\left(\bigcup_{j=i+1}^{k}T_{K}^{n}(C_{j}),
\bigcup_{j=i+1}^{k}B_{j}\right)<\delta$;
\item $d(S^{n}(y), y_{2})<\delta$.
\end{enumerate}
Therefore,
$$
\left(T_{F}|_{\mathbb{F}^{=\lambda}(X)}\right)^{n}(C)\in \left(T_{F}|_{\mathbb{F}^{=\lambda}(X)}\right)^{n}\left(B_{d_{\infty}}(A, \delta)\cap\mathbb{F}^{=\lambda}(X)\right)
\cap \left[B_{d_{\infty}}(B, \delta)\cap\mathbb{F}^{=\lambda}(X)\right]\neq \emptyset,
$$
and
$$
S^{n}(y)\in S^{n}(B(y_{1}, \delta))\cap B(y_{2}, \delta)\neq \emptyset,
$$
implying that
$$
\left(T_{F}|_{\mathbb{F}^{=\lambda}(X)}\times S\right)^{n}(W)\cap V\neq \emptyset.
$$
Hence, $T_{F}|_{\mathbb{F}^{=\lambda}(X)}\times S$ is transitive.
\epf

\medskip

Summing up Lemma \ref{T=T_K}, Lemma \ref{T_F-->T} and Lemma \ref{T_K-->T_F}, one has

\begin{thm}\label{mild}
Let $(X, T)$ be a dynamical system. Then, the following statements are equivalent:
\begin{enumerate}[(1)]
\item $(X, T)$ is mildly mixing;
\item $(K(X), T_{K})$ is mildly mixing;
\item $\forall \lambda\in (0, 1]$, $(\mathbb{F}^{=\lambda}(X), T_{F}|_{\mathbb{F}^{=\lambda}(X)})$ is mildly mixing;
\item $\exists \lambda\in (0, 1]$, $(\mathbb{F}^{=\lambda}(X), T_{F}|_{\mathbb{F}^{=\lambda}(X)})$ is mildly mixing.
\end{enumerate}
\end{thm}

Applying Lemma \ref{a-tran} and Lemma \ref{K(X)},
similarly to the proof of Lemma \ref{T_K-->T_F}, it can be verified
that the following holds.

\begin{thm}\label{mild}
Let $(X, T)$ be a dynamical system. Then, the following statements are equivalent:
\begin{enumerate}[(1)]
\item $(X, T)$ weakly mixing and $\mathbf{a}$-transitive;
\item $(K(X), T_{K})$ is $\mathbf{a}$-transitive;
\item $\forall \lambda\in (0, 1]$, $(\mathbb{F}^{=\lambda}(X), T_{F}|_{\mathbb{F}^{=\lambda}(X)})$ is $\mathbf{a}$-transitive;
\item $\exists \lambda\in (0, 1]$, $(\mathbb{F}^{=\lambda}(X), T_{F}|_{\mathbb{F}^{=\lambda}(X)})$ is $\mathbf{a}$-transitive.
\end{enumerate}
\end{thm}

\section{Equicontinuity of $(\mathbb{F}_{0}(X), T_{F})$}\label{S-5}

Based on Bauer and Sigmund's result which states that $(X, T)$ is equicontinuous
if and only if $(K(X), T_{K})$ is equicontinuous, this section proves the equivalence of
the equicontinuity between $(X, T)$ and $(\mathbb{F}_{0}(X), T_{F})$.

\begin{lem}\cite[Proposition 7]{BS75}\label{X=K(X)}
A dynamical system $(X, T)$ is equicontinuous if and only if $(K(X), T_{K})$ is equicontinuous.
\end{lem}

\begin{prn}\label{dense-equi}
Let $(X, d)$ be a complete metric space, $T: X\longrightarrow X$ be a continuous map and $A\subset X$
be a dense subset of $X$. If $T|_{A}: A\longrightarrow X$ is equicontinuous, then $T: X\longrightarrow X$
is equicontinuous.
\end{prn}

\pf
For any fixed $\varepsilon>0$, there exists $\delta>0$ such that for any $x, y\in A$ with $d(x, y)<\delta$
and any $n\in \mathbb{Z}^+$,
\begin{equation}\label{e-1}
d(T^{n}(x), T^{n}(y))<\frac{\varepsilon}{2}.
\end{equation}
For any $x', y'\in X$ with $d(x', y')<\frac{\delta}{2}$ and any $n\in \mathbb{Z}^{+}$,
as $T^{n}$ is continuous at $x'$ and $y'$, there exists $0<\delta'<\frac{\delta}{4}$
such that for any $z_{1}\in B(x', \delta')$ and any $z_{2}\in B(y', \delta')$,
$$
d(T^{n}(z_{1}), T^{n}(x'))<\frac{\varepsilon}{4} \text{ and } d(T^{n}(z_{2}), T^{n}(y'))<\frac{\varepsilon}{4}.
$$
Choose $p\in B(x', \delta')\cap A$ and $q\in B(y', \delta')\cap A$. Clearly,
$$
d(p, q)\leq
d(p, x')+ d(x', y')+d(y', q)<\delta'+\frac{\delta}{2}+\delta'<\delta.
$$
This, together with
(\ref{e-1}), implies that
$$
d(T^{n}(p), T^{n}(q))<\frac{\varepsilon}{2}.
$$
Then,
\begin{eqnarray*}
d(T^{n}(x'), T^{n}(y'))&\leq & d(T^{n}(x'), T^{n}(p))+d(T^{n}(p), T^{n}(q))+d(T^{n}(q), T^{n}(y'))\\
&<& \frac{\varepsilon}{4}+\frac{\varepsilon}{2}+\frac{\varepsilon}{4}=\varepsilon.
\end{eqnarray*}
Hence, $T$ is equicontinuous.
\epf

\begin{lem}\label{K(X)-->F(X)}
Let $(X, T)$ be a dynamical system. If $(K(X), T_{K})$ is equicontinuous, then $(\mathbb{F}_{0}(X), T_{F})$
is equicontinuous.
\end{lem}

\pf
For any fixed $\varepsilon>0$, as $T_{K}$ is equicontinuous, there exists
$0<\delta<\frac{\mathrm{diam}(X)}{2}$ such that for any $E, F\in K(X)$ with $d_{H}(E, F)<\delta$
and any $n\in \mathbb{Z}^+$, $d_{H}(T_{K}^{n}(E), T_{K}^{n}(F))<\frac{\varepsilon}{2}$. For any
$A, B\in \mathbb{F}_{0}(X)$ with $d_{\infty}(A, B)<\delta$, it is clear that
$\beta:=\max\{A(x): x\in X\}=\max\{B(x): x\in X\}$, as $\delta<\frac{\mathrm{diam}(X)}{2}$,
implying that for any $\alpha\in (0, \beta]$, $[A]_{\alpha} \neq \emptyset \neq [B]_{\alpha}$
and $d_{H}([A]_{\alpha}, [B]_{\alpha})<\delta$.
Then, for any $n\in \mathbb{Z}^+$,
\begin{eqnarray*}
d_{\infty}(T_{F}^{n}(A), T_{F}^{n}(B))&=&\sup\left\{d_{H}(T_{K}^{n}([A]_{\alpha}), T_{K}^{n}([B]_{\alpha})):
\alpha\in (0, 1]\right\}\\
&=&\sup\left\{d_{H}(T_{K}^{n}([A]_{\alpha}), T_{K}^{n}([B]_{\alpha})):
\alpha\in (0, \beta]\right\}
\leq \frac{\varepsilon}{2}<\varepsilon.
\end{eqnarray*}
So, $(\mathbb{F}_{0}(X), T_{F})$
is equicontinuous.
\epf

\begin{lem}\label{F(X)-->K(X)}
If $(\mathbb{F}_{0}(X), T_{F})$ is equicontinuous, then
$(K(X), T_{K})$ is equicontinuous.
\end{lem}

\pf
The result yields by \eqref{e-4}.
\epf

\begin{thm}
The following statements are equivalent:
\begin{enumerate}[(1)]
\item $(X, T)$ is equicontinuous;
\item $(K(X), T_{K})$ is equicontinuous;
\item $(\mathbb{F}_{0}(X), T_{F})$
is equicontinuous.
\end{enumerate}
\end{thm}

\pf
Applying Lemma \ref{X=K(X)}, Lemma \ref{K(X)-->F(X)} and Lemma \ref{F(X)-->K(X)}, this holds trivially.
\epf

\section{Uniform rigidity and proximality of $(\mathbb{F}_{0}(X), T_{F})$}\label{S-6}

Let $n\in \mathbb{N}$, according to Glasner and Maon \cite{GM},
 a dynamical $(X, T)$ is
\begin{enumerate}[(1)]
\item {\it $n$-rigid} if every $n$-tuple $(x_{1}, x_{2}, \ldots, x_{n})\in X^{n}$ is a recurrent
point of $T^{(n)}$;
\item {\it weakly rigid} if $(X, T)$ is $n$-rigid for any $n\in \mathbb{N}$;
\item {\it rigid} if there exists $m_i\longrightarrow +\infty$
such that $T^{m_i} \longrightarrow \mathrm{id}_{X}$ pointwise, where $\mathrm{id}_{X}$ is the
identity map on $X$;
\item {\it uniformly rigid} if there exists $m_{i}\longrightarrow +\infty$
such that $T^{m_{i}}\longrightarrow \mathrm{id}_{X}$ uniformly on $X$.
\end{enumerate}

It can be verified that a dynamical system
$(X, T)$ is uniformly rigid if and only if for any $\varepsilon>0$, there exists $n\in \mathbb{N}$ such that for any
$x\in X$, $d(T^{n}(x), x)<\varepsilon$. It is known that every transitive map containing an
equicontinuous point\footnote{A point $x\in X$ is {\it equicontinuous} if for
any $\varepsilon>0$, there exists $\delta>0$ such that for any $y\in X$ with $d(x, y)<\delta$
and any $n\in \mathbb{Z}^+$, $d(T^{n}(x), T^{n}(y))<\varepsilon$.} is uniformly rigid.
Recently, we \cite{Wu} proved that a dynamical system $(X, T)$
is uniformly rigid if and only if $(M(X), T_{M})$ is uniformly rigid.
The following result is obtained by Li et al. \cite{LOYZ}.

\begin{lem}\cite[Theorem 4.2]{LOYZ}\label{LOYZ}
Let $(X, T)$ be a dynamical system. Then, the following statements are equivalent:
\begin{enumerate}[(1)]
\item $(X, T)$ is uniformly rigid;
\item $(K(X), T_{K})$ is uniformly rigid;
\item $(K(X), T_{K})$ is rigid;
\item $(K(X), T_{K})$ is weakly rigid.
\end{enumerate}
\end{lem}

\begin{prn}\label{Lemma 14}
Let $(X, d)$ be a complete metric space, $T: X\longrightarrow X$ be a continuous map and $A\subset X$
be an invariant dense subset of $X$. If $T|_{A}: A\longrightarrow A$ is uniformly rigid, then $T: X\longrightarrow X$
is uniformly rigid.
\end{prn}

\pf
For any $\varepsilon>0$, as $T|_{A}$ is uniformly rigid, then there exists $n\in \mathbb{N}$ such that for any
$x\in A$, $d(T^{n}(x), x)<\frac{\varepsilon}{2}$. For any $y\in X$, as $T^{n}$ is continuous and $\overline{A}=X$,
there exists $x'\in A$ such that $d(x', y)<\frac{\varepsilon}{4}$
and $d(T^{n}(y), T^{n}(x'))<\frac{\varepsilon}{4}$. This implies that
$$
d(T^{n}(y), y)\leq d(T^{n}(y), T^{n}(x'))+ d(T^{n}(x'), x')+d(x', y)<\frac{\varepsilon}{4}+\frac{\varepsilon}{2}+
\frac{\varepsilon}{4}=\varepsilon.
$$
Thus, $(X, T)$ is uniformly rigid.
\epf


\begin{lem}\label{Lem-K(X)-->F(X)}
If $(K(X), T_{K})$ is uniformly rigid, then $(\mathbb{F}_{0}(X), T_{F})$ is uniformly rigid.
\end{lem}

\pf
For any $\varepsilon>0$, as $T_{K}$ is uniformly rigid, there exists $n\in \mathbb{N}$ such that
for any $E\in K(X)$, $d_{H}(T_{K}^{n}(E), E)<\frac{\varepsilon}{2}$. For any
$A\in \mathbb{F}_{0}(X)$, 
it can be verified that
$$
d_{\infty}(T_{F}^{n}(A), A)=\sup\left\{d_{H}(T_{K}^{n}([A]_{\alpha}), [A]_{\alpha}): \alpha\in (0, 1]\right\}
\leq \frac{\varepsilon}{2}<\varepsilon.
$$
Therefore, $T_{F}$ is uniformly rigid.
\epf

\begin{thm}
Let $(X, T)$ be a dynamical system. Then, the following statements are equivalent:
\begin{enumerate}[(1)]
\item $(X, T)$ is uniformly rigid;
\item $(K(X), T_{K})$ is uniformly rigid;
\item $(\mathbb{F}_{0}(X), T_{F})$ is uniformly rigid;
\item $\forall \lambda\in (0, 1]$, $(\mathbb{F}^{=\lambda}(X), T_{F}|_{\mathbb{F}^{=\lambda}(X)})$ is uniformly rigid;
\item $\forall \lambda\in (0, 1]$, $(\mathbb{F}^{\geq \lambda}(X), T_{F}|_{\mathbb{F}^{=\lambda}(X)})$ is uniformly rigid;
\item $\exists \lambda\in (0, 1]$, $(\mathbb{F}^{=\lambda}(X), T_{F}|_{\mathbb{F}^{=\lambda}(X)})$ is uniformly rigid;
\item $\exists \lambda\in (0, 1]$, $(\mathbb{F}^{\geq \lambda}(X), T_{F}|_{\mathbb{F}^{=\lambda}(X)})$ is uniformly rigid.
\end{enumerate}
\end{thm}

\pf
This follows by Lemma \ref{LOYZ}, Lemma \ref{Lem-K(X)-->F(X)} and \eqref{e-4}.
\epf

\medskip

Recall that a dynamical system $(X, T)$ is proximal if for any $x, y\in X$, $\liminf_{n\rightarrow\infty}
d(T^{n}(x), T^{n}(y))=0$.

\begin{thm}
Let $(X, T)$ be a dynamical system. Then, the following statements are equivalent:
\begin{enumerate}[(1)]
\item $(K(X), T_{K})$ is proximal;
\item $\lim_{n\rightarrow\infty}\mathrm{diam}(T^{n}(X))=0$;
\item $\forall \lambda\in (0, 1]$, $(\mathbb{F}^{=\lambda}(X), T_{F}|_{\mathbb{F}^{=\lambda}(X)})$ is proximal.
\end{enumerate}
\end{thm}

\pf
(3) $\Longrightarrow$ (1) follows by Remark \ref{Remark 3}.

(1) $\Longrightarrow$ (2). Fix any $x\in X$. As $T_{K}$ is proximal, $\liminf_{n\rightarrow \infty}
d_{H}\left(T_{K}^{n}(X), T_{K}^{n}(\{x\})\right)=0$. This, together with the fact that $\{T^n (X)\}$ is a
decreasing sequence, implies that $\lim_{n\rightarrow\infty}\mathrm{diam}(T^{n}(X))=0$.

(2) $\Longrightarrow$ (3). For nay $A, B\in \mathbb{F}^{=\lambda}(X)$ and any $n\in \mathbb{Z}^{+}$,
one has
\begin{eqnarray*}
d_{\infty}\left(\left(T_{F}|_{\mathbb{F}^{=\lambda}(X)}\right)^{n}(A),
\left(T_{F}|_{\mathbb{F}^{=\lambda}(X)}\right)^{n}(B)\right)
&=&\sup\left\{d_{H}(T^{n}([A]_{\alpha}), T^{n}([B]_{\alpha})): \alpha\in (0, 1]\right\}\\
&=&\sup\left\{d_{H}(T^{n}([A]_{\alpha}), T^{n}([B]_{\alpha})): \alpha\in (0, \lambda]\right\}
\leq \mathrm{diam}(T^{n}(X)).
\end{eqnarray*}
This implies that $\liminf_{n\rightarrow \infty}d_{\infty}\left(\left(T_{F}|_{\mathbb{F}^{=\lambda}(X)}\right)^{n}(A),
\left(T_{F}|_{\mathbb{F}^{=\lambda}(X)}\right)^{n}(B)\right)=0$.
\epf

\begin{rmk}
\begin{enumerate}[(1)]
\item Because there exists proximal system which is topologically mixing, this shows that the proximality
of $(\mathbb{F}^{=\lambda}(X), T_{F}|_{\mathbb{F}^{=\lambda}(X)})$ is strictly stronger that the proximality of
$(X, T)$.

\item Clearly, for any $A, B\in \mathbb{F}_{0}(X)$ with $\xi:=\max\{A(x): x\in X\}<
\max\{B(x): x\in X\}:=\eta$ and any $n\in \mathbb{Z}^{+}$,
$$
d_{\infty}(T_{F}^{n}(A), T_{F}^{n}(B))\geq d_{H}(T^{n}([A]_{\frac{\xi+\eta}{2}}),
T^{n}([B]_{\frac{\xi+\eta}{2}}))=d_{H}(\emptyset, T^{n}([B]_{\frac{\xi+\eta}{2}}))
=\mathrm{diam}(X).
$$
This implies that $(\mathbb{F}_{0}(X), T_{F})$ is not proximal.
\end{enumerate}
\end{rmk}

\section*{Acknowledgments}
We thank the referees for their careful reading and
valuable suggestions which helped us to improve the quality of this paper.

\end{document}